\documentclass[runningheads,a4paper]{llncs}

\usepackage[dvips]{epsfig}
\usepackage{amsfonts}
\usepackage{amssymb}
\usepackage{amscd}
\usepackage{amstext}
\usepackage{amsmath}
\usepackage{pstricks}
\usepackage{enumerate}

\newtheorem{corollaire}[example]{Corollary}

\def\Proof{\medskip\noindent {\it Proof --- \ }}
\def\cqfd{\hfill $\Box$ \bigskip}

\def\hs{\hbox to 3mm{}}
\def\hhs{\hbox to 5cm{}}
\def\ssk{\smallskip}

\def\shuff#1#2{\mathbin{
      \hbox{\vbox{
        \hbox{\vrule
              \hskip#2
              \vrule height#1 width 0pt
               }%
        \hrule}%
             \vbox{
        \hbox{\vrule
              \hskip#2
              \vrule height#1 width 0pt
               \vrule }%
        \hrule}%
}}}
\def\shuffle{{\mathchoice{\shuff{7pt}{3.5pt}}%
                        {\shuff{6pt}{3pt}}%
                        {\shuff{4pt}{2pt}}%
                        {\shuff{3pt}{1.5pt}}}}

\def\adots{\mathinner{\mkern2mu\raise1pt\hbox{.}
\mkern3mu\raise4pt\hbox{.}\mkern1mu\raise7pt\hbox{.}}}

\def\up#1{\raise 1ex\hbox{\footnotesize#1}}

\def\mref#1{(\ref{#1})}

\def\H{\mathcal{H}}
\def\M{\mathcal{M}}
\def\bd{\mathbf{d}}

\def\supp{\mathop\mathrm{supp}\nolimits}
\def\span{\mathop\mathrm{span}\nolimits}

\def\ra{\rightarrow}

\def\path{\rightsquigarrow}

\def\A{\mathcal{A}}

\def\B{\mathcal{B}}

\def\F{\mathcal{F}}
\def\N{{\mathbb N}}
\def\C{{\mathbb C}}
\def\R{{\mathbb R}}
\def\Z{{\mathbb Z}}

\def\X{{\mathbb X}}
\def\al{\alpha}
\def\be{\beta}

\def\calC{{\mathcal C}}
\def\scal#1#2{\langle #1 | #2 \rangle}

\def\ncs#1#2{#1\langle\langle #2\rangle\rangle}

\def\ncp#1#2{#1\langle #2\rangle}

\def\End{\mathrm{End}}
\def\Aut{\mathrm{Aut}}

\def\lp#1#2{\mathrm{Lie}_#1\langle #2\rangle}
\def\ls#1#2{\mathrm{Lie}_#1\langle\langle #2\rangle\rangle}
\def\Li{\mathrm{Li}}

\def\calA{\mathcal{A}}

\usepackage{color}
\usepackage{ulem}

\begingroup
\def\2#1{\ifnum#1<10 0\fi\the#1}
\xdef\isodayandtime{
{\2\day-\2\month-\the\year\space\2{\count0}:%
\2{\count2}}}
\endgroup

\title{\Large\bf Independence of hyperlogarithms
over function fields via algebraic combinatorics.}
\titlerunning{Independence of hyperlogarithms}
\author{
{\sc M. Deneufch\^atel, G. H. E. Duchamp, Hoang Ngoc Minh }\\
{\rm \small \it Institut Galil\'ee, LIPN - UMR 7030
CNRS - Universit\'e Paris 13
F-93430 Villetaneuse, France}\\
{\rm \small and} \\
{\sc  Allan I. Solomon }\\
{\rm \small \it Department of Physics and Astronomy, The Open University, Milton Keynes MK7 6AA, UK\\
and\\
 LPTMC - UMR 7600 CNRS - Universit\'e Paris 6 
F-75252 Paris, France}
}
\authorrunning{M. Deneufch\^atel, G. H. E. Duchamp, Hoang Ngoc Minh 
\and Allan I. Solomon}

\institute{Institut Galil\'ee \and Open University}
\date{}

\begin{document}

\maketitle


\begin{abstract}
We obtain a  necessary and sufficient condition for the linear independence of solutions of differential equations for  hyperlogarithms. The key fact is that the multiplier
(i.e. the factor $M$ in the differential equation $dS=MS$) has only singularities of first order (Fuchsian-type equations) and this implies that they freely span a space which contains no primitive. We give direct applications where we extend the property of linear independence to the largest known ring of coefficients.

\end{abstract}

\section{Introduction}

In his 1928 study of the solutions of linear differential equations following Poincar\'e, Lappo-Danilevski introduced the so-called {\it hyperlogarithmic
functions of order $m$}, functions of iterated integrals of the following form with logarithmic poles \cite{lappo}~:
\begin{eqnarray}
L(a_0,\ldots,a_n|z,z_0)=\int_{z_0}^z\int_{z_0}^{s_n}\ldots\int_{z_0}^{s_1}
\frac{ds_0}{s_0-a_0}\ldots\frac{ds_n}{s_n-a_n},
\end{eqnarray}
where $z_0$ is a fixed point. It suffices that $z_0 \neq a_0$ for this iterated integral to converge.
The classical polylogarithm $\Li_n$ is a particular case of these integrals \cite{lewin1}~:
\begin{eqnarray}
\Li_n(z)=\int_0^z\int_0^{s_n}\ldots\int_0^{s_2}
\frac{ds_1}{1-s_1}\frac{ds_2}{s_2}\ldots\frac{ds_n}{s_n}
=-L(1,\underbrace{0,\ldots,0}_{n-1\ \mathrm{times}}|z, 0).
\end{eqnarray}
These iterated integrals also appear in quantum electrodynamics (see \cite{dyson,magnus} for example).
Chen \cite{chen} studied them systematically and provided a  noncommutative algebraic context in which to treat them.
Fliess \cite{fliess1,fliess2} encoded these iterated integrals by words over a finite alphabet
and extended them to a symbolic calculus\footnote{A kind of Feynman like
operator calculus \cite{feynman}.} for nonlinear differential equations
of the following form, in the context of noncommutative formal power series~:
\begin{eqnarray}\label{nonlinear}
\left\{\begin{array}{lcl}
y(z)&=&f(q(z)),\\
\dot q(z)&=&{\displaystyle\sum_{i=0}^m\frac{A_i(q)}{z-a_i}},\\
q(z_0)&=&q_0,
\end{array}\right.
\end{eqnarray}
where the state $q=(q_1,\ldots,q_n)$ belongs to a complex analytic manifold of dimension $N$,
$q_0$ denotes the initial state, the observable $f$ belongs to $\C^{\rm cv}[\![q_1,\ldots,q_N]\!]$,
and $\{A_i\}_{i=0,n}$ is the polysystem defined as follows
\begin{eqnarray}\label{vectorfield}
A_i(q)=\sum_{j=1}^nA_i^j(q)\frac{\partial}{\partial q_j},
\end{eqnarray}
with, for any $j=1,\ldots,n$,
$A_i^j(q)\in\C^{\rm cv}[\![q_1,\ldots,q_N]\!]$.\\
By introducing the encoding alphabet $X=\{x_0,\ldots,x_m\}$,
the method of Fliess consists in exhibiting two formal power series over the monoid $X^*$~:
\begin{eqnarray}\label{chenfliess}
F:=\sum_{w\in X^*}\calA(w)\circ f_{|_{q_0}}\;w&\mbox{and}&
C:=\sum_{w\in X^*}\alpha_{z_0}^z(w)\;w
\end{eqnarray}
in order to compute the output $y$. These series are subject to convergence conditions (precisely speaking, the convergence of a duality pairing), as follows:
\begin{eqnarray}\label{output}
y(z)=\langle F||C\rangle:=\sum_{w\in X^*}\calA(w)\circ f_{|_{q_0}}\;\alpha_{z_0}^z(w),
\end{eqnarray}
where
\begin{itemize}
\item $\calA$ is a morphism of algebras from $\ncs{\C}{X}$ to the algebra generated by the polysystem $\{A_i\}_{i=0,n}$~:
\begin{eqnarray}
\calA(1_{X^*})&=&\mbox{identity},\\
\forall w=vx_i,x_i\in X,v\in X^*,\quad
\calA(w)&=&\calA(v)A_i
\end{eqnarray}
\item $\alpha_{z_0}^z$ is a shuffle algebra morphism from $(\ncs{\C}{X},\shuffle)$ to some differential field $\calC$~:
\begin{eqnarray}
\alpha_{z_0}^z(1_{X^*})&=&1,\\
\forall w=vx_i,x_i\in X,v\in X^*,\quad
\alpha_{z_0}^z(w)&=&\int_{z_0}^z\frac{\alpha_{z_0}^s(v)}{s-a_i}.
\end{eqnarray}
\end{itemize}
Formula (\ref{output}) also states  that the iterated integrals over the rational functions
\begin{eqnarray}
u_i(z)=\frac1{z-a_i},&&i=0,..,n,
\end{eqnarray}
span the vector space $\calC$.

As for the linear differential equations, the essential difficulty is to construct the fundamental system of solutions, or the Picard-Vessiot extension, to describe the space of solutions of the differential system (\ref{nonlinear}) algorithmically  \cite{singervanderput}. For that, one needs to prove the linear independence of the iterated integrals in order to obtain the {\it universal} Picard-Vessiot extension.
The $\C$-linear independence was already  shown by Wechsung \cite{wechsung}. His method uses a recurrence based on the total degree. However this method cannot be used with
variable coefficients. Another proof was given in
\cite{FPSAC98} based on monodromy.
In this note  we describe  a general theorem in differential computational algebra
and show that, at the cost of using variable domains (which
is the realm of germ spaces), and replacing the recurrence on total
degree by a recursion on the words (with graded lexicographic ordering),
one can encompass the previous results mentioned above and obtain much larger
rings of coefficients and configuration alphabets (even infinite of continuum cardinality).

\smallskip
{\bf Acknowledgements} : The authors are pleased to acknowledge the hospitality of institutions in Paris and UK. We take advantage of these lines to acknowledge support from the French Ministry of Science and Higher Education under Grant ANR PhysComb and local support from the Project "Polyzetas".

\section{Non commutative differential equations.}

We recall the Dirac-Sch\"utzenberger notation, as in \cite{BR,D21,Reu}. Let $X$ be an alphabet and $R$ be a commutative ring with unit. The algebra of noncommutative polynomials is the algebra $R[X^*]$ of the free monoid $X^*$. As an $R$-module, $R^{(X^*)}$ is the set of finitely supported $R$-valued function on $X^*$ and, as such, it is in natural duality with the algebra of all functions on $X^*$ (the large algebra of $X^*$ \cite{B_Alg_III}), $R^{X^*}=\ncs{R}{X}$, the duality being given, for $f\in \ncs{R}{X}$ and $g\in R[X^*]$, by
\begin{equation}
	\scal{f}{g}=\sum_{w\in X^*}f(w)g(w)\ .
\end{equation}

The r\^ole of the ring is played here by a commutative differential $k$-algebra $(\A,d)$; that is, a $k$-algebra $\A$ (associative and commutative with unit)  endowed with a distinguished derivation $d\in \frak{Der}(\A)$ (the ground field $k$ is supposed commutative and of characteristic zero). We assume that the ring of constants $ker(d)$ is precisely $k$.

\ssk
An alphabet $X$ being given, one can at once extend the derivation $d$ to a derivation of the algebra $\ncs{\A}{X}$ by
\begin{equation}
	\bd(S)=\sum_{w\in X^*}d(\scal{S}{w})w\ .
\end{equation}
We now act with  this derivation $\bd$ on the power series $C$ given in (\ref{chenfliess}). We then get~:
\begin{equation}
	\bd(C)=\biggl(\sum_{i=1}^mu_i x_i\biggr)C\ .
\end{equation}

We are now in a position to state the main theorem which resolves many important questions, some of which we  shall see  in the applications.

\begin{theorem}\label{ind_lin} Let $(\A,d)$ be a $k$-commutative associative differential algebra with unit ($ch(k)=0, ker(d)=k$) and $\calC$ be a differential subfield of $\A$ (i.e. $d(\calC)\subset \calC$).
We suppose that $S\in \ncs{\A}{X}$ is a solution of the differential equation
\begin{equation}
	\bd(S)=MS\ ;\ \scal{S}{1}=1
\end{equation}
where the multiplier $M$ is a homogeneous series (a polynomial in the case of finite $X$) of degree $1$,
i.e.
\begin{equation}
M=\sum_{x\in X}u_x x\in \ncs{\calC}{X}\ .
\end{equation}
The following conditions are equivalent :
\begin{enumerate}[i)]
	\item The family $(\scal{S}{w})_{w\in X^*}$ of coefficients of $S$ is free over $\calC$.
  \item The family of coefficients $(\scal{S}{y})_{y\in X\cup \{1_{X^*}\}}$ is free over $\calC$.
	\item The family $(u_x)_{x\in X}$ is such that, for $f\in \calC$ and $\al_x\in k$
\begin{equation}\label{prim_ind}
d(f)=\sum_{x\in X} \al_x u_x \Longrightarrow (\forall x\in X)(\al_x=0)\ .
\end{equation}
	\item The family $(u_x)_{x\in X}$ is free over $k$ and
\begin{equation}
	d(\calC)\cap\span_k\Big((u_x)_{x\in X}\Big)=\{0\}\ .
\end{equation}
\end{enumerate}
\end{theorem}

\Proof (i)$\Longrightarrow$(ii) Obvious.\\
(ii)$\Longrightarrow$(iii)\\
Suppose that the family $(\scal{S}{y})_{y\in X\cup \{1_{X^*}\}}$ (coefficients taken at letters and the empty word) of coefficients of $S$ is free over $\calC$ and let us consider the relation as in eq. \mref{prim_ind}
\begin{equation}
d(f)=\sum_{x\in X} \al_x u_x\ .
\end{equation}
We form the polynomial $P=-f1_{X^*}+\sum_{x\in X} \al_x x$. One has $\bd(P)=-d(f)1_{X^*}$ and
\begin{equation}
	d(\scal{S}{P})=\scal{\bd(S)}{P}+\scal{S}{\bd(P)}=\scal{MS}{P}-d(f)\scal{S}{1_{X^*}}=
	(\sum_{x\in X} \al_x u_x)-d(f)=0
\end{equation}
whence  $\scal{S}{P}$ must be a constant, say $\lambda\in k$. For $Q=P-\lambda.1_{X^*}$, we have
$$\supp(Q)\subset X\cup \{1_{X^*}\}\ \textrm{and}\ \scal{S}{Q}=\scal{S}{P}-\lambda\scal{S}{1_{X^*}}=\scal{S}{P}-\lambda=0
\ .
$$
This implies that $Q=0$ and, as $Q=-(f+\lambda)1_{X^*}+\sum_{x\in X} \al_x x$, one has, in particular, all the $\al_x=0$.\\
(iii)$\Longleftrightarrow$(iv)\\
Obvious, (iv) being a geometric reformulation of (iii).\\
(iii)$\Longleftrightarrow$(i)\\
Let $\mathcal{K}$ be the kernel of $P\mapsto \scal{S}{P}$ (a linear form
$\ncp{\calC}{X}\ra \calC$) i.e.
\begin{equation}
	\mathcal{K}=\{P\in \ncp{\calC}{X}| \scal{S}{P}=0\}\ .
\end{equation}
If $\mathcal{K}=\{0\}$, we are done. Otherwise, let us adopt the following strategy.\\
First, we order $X$ by some well-ordering $<$ (\cite{B_Set} III.2.1) and $X^*$ by the graded lexicographic ordering
$\prec$ defined by
\begin{equation}
	u\prec v \Longleftrightarrow |u|<|v|\ \textrm{or}\ (u=pxs_1\ ,\ v=pys_2\ \textrm{and}\ x<y).
\end{equation}
It is easy to check that $\prec$ is also a well-ordering relation. For each nonzero polynomial $P$, we denote by $lead(P)$ its leading monomial; i.e. the greatest element of its support $\supp(P)$ (for $\prec$).\\
Now, as $\mathcal{R}=\mathcal{K}-\{0\}$ is not empty, let $w_0$ be the minimal element of $lead(\mathcal{R})$ and choose a $P\in \mathcal{R}$ such that $lead(P)=w_0$. We write
\begin{equation}
	P=fw_0+\sum_{u\prec w_0}\scal{P}{u}u\ ;\ f\in \calC-\{0\}\ .
\end{equation}
The polynomial $Q=\frac{1}{f}P$ is also in $\mathcal{R}$ with the same leading monomial, but the leading coefficient is now $1$; and so $Q$ is given by
\begin{equation}
	Q=w_0+\sum_{u\prec w_0}\scal{Q}{u}u\ .	
\end{equation}

Differentiating $\scal{S}{Q}=0$, one gets
\begin{eqnarray}
	0=\scal{\bd(S)}{Q}+\scal{S}{\bd(Q)}=\scal{MS}{Q}+\scal{S}{\bd(Q)}=
	\scal{S}{M^\dagger Q}+\scal{S}{\bd(Q)}=\scal{S}{M^\dagger Q+\bd(Q)}
\end{eqnarray}
with
\begin{equation}
M^\dagger Q+\bd(Q)=\sum_{x\in X} u_x (x^\dagger Q)+\sum_{u\prec w_0}d(\scal{Q}{u})u\in \ncp{\calC}{X}	\ .
\end{equation}
It is impossible that $M^\dagger Q+\bd(Q)\in \mathcal{R}$ because it would be of leading monomial strictly less than $w_0$, hence $M^\dagger Q+\bd(Q)=0$. This is equivalent to the recursion
\begin{equation}
	d(\scal{Q}{u})=-\sum_{x\in X} u_x \scal{Q}{xu}\ ;\ \textrm{for}\ x\in X\ ,\ v\in X^*   .
\end{equation}
From this last relation, we deduce that $\scal{Q}{w}\in k$ for every $w$ of length $deg(Q)$ and,
because $\scal{S}{1}=1$, one must have $deg(Q)>0$. Then, we write $w_0=x_0v$ and compute the coefficient at $v$
\begin{equation}
d(\scal{Q}{v})=-\sum_{x\in X} u_x \scal{Q}{xv}=\sum_{x\in X} \al_x u_x
\end{equation}
with coefficients $\al_x=-\scal{Q}{xv}\in k$ as $|xv|=\deg(Q)$ for all $x\in X$.
Condition \mref{prim_ind} implies that all coefficients $\scal{Q}{xu}$ are zero; in particular, as
$\scal{Q}{x_0u}=1$, we get a contradiction. This proves that $\mathcal{K}=\{0\}$.\\
\cqfd

\section{Applications}

Let $V$ be a connected and simply connected analytic variety (for example, the doubly cleft plane\\ $\C-(]-\infty, 0[\cup ]1,+\infty[)$, the Riemann sphere or the universal covering of $\C-\{0,1\}$), and let  $\H=C^\omega(V,\C)$ be the space of analytic functions on $V$.

It is possible to enlarge the range of scalars to coefficients that are analytic functions with variable domains $f:\ dom(f)\ra\C$.

\begin{definition} We define a \textrm{differential field of germs} as the data of a filter basis $\B$ of open connected subsets of $V$, and a map $\calC$ defined on $\B$ such that for every $U\in \B$, $\calC[U]$ is a subring of $C^\omega(U,\C)$ and

\begin{enumerate}
	\item $\calC$ is compatible with restrictions i.e. if $U,V\in\B$ and  $V\subset U$, one has
	$$
	res_{VU}(\calC[U])\subset \calC[V]
	$$
	\item if $f\in \calC[U]\setminus \{0\}$ then there exists $V\in\B$ s.t.
	$V\subset U-\mathcal{O}_f$ and $f^{-1}$ (defined on $V$) is in $\calC[V]$\ .
\end{enumerate}
\end{definition}

There are important cases where the conditions \mref{lin_indep} are satisfied as shown by the following theorem.

\begin{theorem}\label{lin_indep} Let $V$ be a simply connected non-void open subset of $\C-\{a_1,\cdots a_n\}$ ($\{a_1,\cdots a_n\}$ are distinct points), $M=\sum_{i=1}^n \frac{\lambda_i x_i}{z-a_i}$ be a multiplier on $X=\{x_1,\cdots x_n\}$ with all $\lambda_i\not=0$ and $S$ be any regular solution of
\begin{equation}\label{diff_eq2}
	\frac{d}{dz}S=MS\ .
\end{equation}
Then, let $\calC$ be a differential field of functions defined on $V$ which does not contain linear combinations of logarithms on any domain but which contains $z$ and the constants (as, for example the rational functions).\\
If $U\in \B$ (i.e. $U$ is a domain of $\calC$) and $P\in \ncp{\calC[U]}{X}$, one has
\begin{equation}
	\scal{S}{P}=0\Longrightarrow P=0
\end{equation}
\end{theorem}

\Proof Let $U \in\B$. For every non-zero $Q\in \ncp{\calC[U]}{X}$, we denote by $lead(Q)$ the greatest word in the support of $Q$ for the graded lexicographic ordering $\prec$.  We endow $X$ with an arbitrary linear ordering,  and call $Q$ monic if the leading coefficient $\scal{Q}{lead(Q)}$ is $1$. A monic polynomial is then given by
\begin{equation}
	Q=w+\sum_{u\prec w}\scal{Q}{u}u\ .
\end{equation}
Now suppose  that it is possible to find $U$ and $P\in\ncp{\calC[U]}{X}$ (not necessarily monic) such that $\scal{S}{P}=0$;  we choose $P$ with $lead(P)$ minimal for $\prec$.

\ssk
Then
\begin{equation}
	P=f(z)w+\sum_{u\prec w}\scal{P}{u}u
\end{equation}
with $f\not\equiv 0$. Thus $U_1=U\setminus\mathcal{O}_f\in \B$
and $Q=\frac{1}{f(z)}P\in \ncp{\calC[U_1]}{X}$ is monic and satisfies
\begin{equation}\label{linear_rel}
\scal{S}{Q}=0\ .	
\end{equation}
Differentiating eq. \mref{linear_rel}, we get
\begin{equation}
0=\scal{S'}{Q}+\scal{S}{Q'}=\scal{MS}{Q}+\scal{S}{Q'}=\scal{S}{Q'+M^\dagger Q}\ .	
\end{equation}
Remark that one has
\begin{equation}
	Q'+M^\dagger Q\in \ncp{\calC[U_1]}{X}
\end{equation}
If $Q'+M^\dagger Q\not= 0$, one has $lead(Q'+M^\dagger Q)\prec lead(Q)$ and this is not
possible because of the minimality hypothesis of $lead(Q)=lead(P)$. Hence, one must have
$R=Q'+M^\dagger Q= 0$. With $|w|=n$, we now write
\begin{equation}
	Q=Q_n+\sum_{|u|< n}\scal{Q}{u}u\
\end{equation}
where $Q_n=\sum_{|u|=n}\scal{Q}{u}u$ is the dominant homogeneous component of $Q$.
For every $|u|=n$ we have
\begin{equation}
	(\scal{Q}{u})'=-\scal{M^\dagger Q}{u}=-\scal{Q}{Mu}=0
\end{equation}
thus all the coefficients of $Q_n$ are constant.

\ssk
If $n=0$, $Q\not=0$ is constant which is impossible by eq. \mref{linear_rel} and because $S$ is regular. If $n>0$, for any word $|v|=n-1$, we have
\begin{equation}
	(\scal{Q}{v})'=-\scal{M^\dagger Q}{v}=-\scal{Q}{Mv}=
	-\sum_{i=0}^n \frac{\lambda_i}{z-a_i} \scal{Q}{x_iv}=
	-\sum_{i=0}^n \frac{\lambda_i}{z-a_i} \scal{Q_n}{x_iv}
\end{equation}
bcause all $x_iv$ are of length $n$.\\
Then
\begin{equation}
	\scal{Q}{v}=-\sum_{i=0}^n \scal{Q_n}{x_iv} \int_{\al}^z \frac{\lambda_i}{s-a_i} ds + const
\end{equation}
But all the functions $\int_{\al}^z \frac{\lambda_i}{s-a_i} ds$ are linearly independent over $\C$ and not all the scalars $\scal{Q_n}{x_iv}$ are zero (write $w=x_kv$ and choose $v$ accordingly).  This contradicts the fact that $Q\in \ncp{\calC[U_1]}{X}$ as $\calC$ contains no linear combination of logarithms.
\cqfd

\begin{corollaire} Let $V$ be as above and $R$ be the ring of functions which can be analytically extended to some $V\cup U_{a_1}\cup U_{a_2}\cup\cdots U_{a_n}$ where $U_{a_i}$ are open neighborhoods of $a_i,i=1\cdots n$ and have non-essential singularities at these points. Then, the set of hyperlogarithms $(\scal{S}{w})_{w\in X^*}$ are linearly independent over $R$.
\end{corollaire}

\begin{remark}
\begin{enumerate}[i)]
	\item If a series $S=\sum_{w\in X^*}\scal{S}{w}w$ is a regular solution of \mref{diff_eq2} and satisfies the equivalent conditions of the theorem
\mref{ind_lin}, then so too does every $Se^C$ (with $C\in\ls{\C}{X}$) .\\
	\item Series such as that  of polylogarithms and all the exponential solutions of equation
\begin{equation}\label{diff_eq3}
	\frac{d}{dz}(S)=(\frac{x_0}{z}+\frac{x_1}{1-z})S
\end{equation}
satisfy the conditions of the theorem \mref{ind_lin} as shown by  theorem \mref{lin_indep}.\\
	\item Call $\F(S)$ the vector space generated by the coefficients of the series $S$. One may ask what happens when the conditions for independence are not satisfied.\\
In fact, the set of Lie series $C\in \ls{\C}{\X}$ such that there exists a $\phi\in\End(\F(S))$ (thus a derivation) such that $SC=\phi(S)$, is a closed Lie subalgebra of $\ls{\C}{\X}$ which we will denote by $Lie_S$. For example
\begin{itemize}
	\item for $X=\{x_0,x_1\}$ and $S=e^{zx_0}$ one has $x_0\in Lie_S\ ;\ x_1\notin Lie_S$
	\item for $X=\{x_0,x_1\}$ and $S=e^{z(x_0+x_1)}$, one has $x_0,x_1\notin Lie_S$ but $(x_0+x_1)\in Lie_S$.
\end{itemize}
	\item Theorem \mref{lin_indep} holds {\it mutatis mutandis} when the multiplier is infinite i.e.
$$
M=\sum_{i\in I} \frac{\lambda_i x_i}{z-a_i}
$$
even if $I$ is continuum infinite (say $I=\R$, singularities being all the reals).
	\item Theorem \mref{lin_indep}  no longer holds with singularities of higher order (i.e. not Fuchsian). For example, with
\begin{equation}	
	M=\frac{x_0}{z^2}+\frac{x_1}{(1-z)^2}\ .
\end{equation}

Firstly, the differential field $\calC$ generated by
\begin{equation}	
u_0=\frac{1}{z^2},\ u_1=\frac{1}{(1-z)^2}
\end{equation}	
contains
\begin{equation}
\frac{d}{dz}(\frac{1}{2u_0})=z
\end{equation}
and hence $\calC=\C(z)$, the field of rational functions over $\C$. Condition (ii) of Theorem
\mref{ind_lin} is not satisfied (as $z^2u_0-(1-z)^2u_1=0$). Moreover, one has also $\Z$-dependent relations such as
\begin{equation}
	\scal{S}{x_1x_0}+\scal{S}{x_0x_1}+\scal{S}{x_1}-\scal{S}{x_0}=0\ .
\end{equation}
\end{enumerate}
\end{remark}

\section{Through the looking glass: passing from right to left.}\label{fonctions_S}

We are still in the context of analytic functions as above. A series $S\in\ncs{\H}{X}$ is said to be group-like if
\begin{equation}
	\Delta(S)=S\otimes S
\end{equation}
where $\Delta$ is the dual of the shuffle product \cite{Reu} defined on series by $\Delta(S)=\sum_{w\in X^*}\scal{S}{w}\Delta(w)$ and on the words by the recursion ($x\in X,u\in X^*$)
\begin{equation}
	\Delta(1_{X^*})=1_{X^*}\otimes 1_{X^*}\ ;\ \Delta(xu)=
	(x\otimes 1_{X^*}+1_{X^*}\otimes x)\Delta(u)
\end{equation}
Let $S\in\ncs{\H}{X}$. We call $\F(S)$ the $\C$-vector space generated by the coefficients of $S$. One has
\begin{equation}
	\F(S)=\{\scal{S}{P}\}_{P\in \ncp{\C}{X}}\ .
\end{equation}

We recall that, for $a\in X$ and $w\in X^*$, the partial degree $|w|_a$ is the 
number of occurrences of $a$ in $w$, it is defined by the recursion 
\begin{equation}
	|1_{X^*}|_a=0\ ;\ |bu|_a=\delta{b,a}+|u|_a\ .
\end{equation}
 
Or course the lenght of the word is the sum of the partial degrees i.e. 
$|w|=\sum_{x\in X}|w|_x$. 
The function $a\mapsto |w|_a$ belongs to $\N^{(X)}$ (finitely supported functions from 
$X$ to $\N$). For $\al\in \N^{(X)}$, we note $\ncp{\C_{\leq \al}}{X}$, the set of 
polynomials $Q\in \ncp{\C}{X}$ such that $supp(Q)\subset X^{\leq \al}$ i.e.
\begin{equation}
	\scal{Q}{w}\not=0\Longrightarrow (\forall x\in X)(|w|_x\leq \al(x))
\end{equation}

In the same way, we consider the filtration by total degree (length) 
\begin{equation}
\ncp{\C_{\leq n}}{X}=\sum_{|\al|\leq n} \ncp{\C_{\leq \al}}{X}\ .
\end{equation}
We  use the following increasing filtrations
\begin{equation}\label{fine_filtr_coeff}
	\F_{\leq \al}(S)=\{\scal{S}{P}\}_{P\in \ncp{\C_{\leq \al}}{X}}\ .
\end{equation}
or
\begin{equation}\label{filtr_coeff}
	\F_{\leq n}(S)=\{\scal{S}{P}\}_{P\in \ncp{\C_{\leq n}}{X}}\ .
\end{equation}

\begin{proposition}\label{kernull} We have the following properties :
\begin{enumerate}[i)]
\item If $T\in\ncs{\C}{X}$ then $\F(ST)\subset \F(S)$ and one has equality if $T$ is invertible.\\
\item If $S$ is group-like, then $\F(S)$ is a unital sub-algebra of $\H$, which is filtered w.r.t. \mref{fine_filtr_coeff} and \mref{filtr_coeff} i.e.
\begin{equation}
\F_{\leq \al}(S)\F_{\leq \be}(S)\subset \F_{\leq \al+\be}(S)
\end{equation}
\end{enumerate}
\end{proposition}

\Proof (i) The space $\F(ST)$ is spanned by the
$$
\scal{ST}{w}=\sum_{uv=w}\scal{S}{u}\scal{T}{v}\in \F(S)
$$
and if $T$ is invertible one has $\F(S)=\F(STT^{-1})\subset \F(ST)$ which proves the equality.\\
ii) If $S$ is group-like, one has
\begin{equation}
	\scal{S}{u}\scal{S}{v}=\scal{S\otimes S}{u\otimes v}=\scal{\Delta(S)}{u\otimes v}=\scal{S}{u\shuffle v}
\end{equation}

In the case when all the functions $\scal{S}{w}$ are $\C$-linearly independent, one has a correspondence between the differential Galois group (acting on the right) of a differential equation of type \mref{diff_eq3} (acting on the right) and the group of automorphisms of $\F(S)$ compatible with the preceding filtration (they turn out to be unipotent).

\newpage
\begin{proposition}\label{ind_lin}
Let $S$ be a group-like series. The following conditions are equivalent:
\begin{enumerate}[i)]
\item For every $x\in X$, $ker_\C(S)\subset ker_\C(Sx)$.\\
\item For every $x\in X$, there is a derivation $\delta_x\in \frak{Der}(\F(S))$ such that
\begin{equation}
	\delta_x(S)=Sx
\end{equation}
\item For every $x\in X$, there is a one-parameter group of automorphisms $\phi_{x}^t\in Aut(\F(S));\ t\in \R$ such that
\begin{equation}
	\phi_{x}^t(S)=Se^{tx}
\end{equation}
\item For every $C\in \ls{C}{X}$, there is $\delta\in \frak{Der}(\F(S))$ such that
\begin{equation}
	\delta(S)=SC
\end{equation}
\item For every $C\in \ls{C}{X}$, there is $\phi\in Aut(\F(S))$ such that
\begin{equation}
	\phi(S)=Se^C
\end{equation}
\item The functions $(\scal{S}{w})_{w\in X^*}$ are $\C$-linearly independant.
\end{enumerate}
\end{proposition}

\Proof $i)\Longrightarrow ii)$ From the inclusion, we deduce that, for all $x\in X$ there exists a $\C$-linear mapping $\phi\in \End(\F(S))$ such that for all $w\in \M$,  $\phi(\scal{S}{w})=\scal{Sx}{w}$. It must be a derivation of $\F(S)$ as
\begin{eqnarray}
\phi(\scal{S}{u}\scal{S}{v})=\phi(\scal{S}{u\shuffle v})=\scal{Sx}{u\shuffle v}=
\scal{S}{(u\shuffle v)x^{-1}}=\cr
\scal{S}{(ux^{-1}\shuffle v)+(u\shuffle vx^{-1})}=\scal{S}{(ux^{-1}\shuffle v)}\scal{S}{(u\shuffle vx^{-1})}=\cr
\scal{Sx}{u}\scal{S}{v}+\scal{S}{u}\scal{Sx}{v}=
\phi(\scal{S}{u})\scal{S}{v}+\scal{S}{u}\phi(\scal{S}{v})
\end{eqnarray}
from the fact that $(\scal{S}{w})_{w\in X^*}$ spans $\F(S)$.\\
$ii)\Longrightarrow iv)$ As $(\scal{S}{w})_{w\in X^*}$ spans $\F(S)$, the derivation $\phi$ is uniquely defined. We denote it by $\delta_x$,  and notice that, in so doing, we have constructed a mapping $\Phi : X\ra \frak{Der}(\F(S))$, which is Lie algebra. Therefore, there is a unique extension of this mapping as a morphism $\lp{\C}{X}\ra \frak{Der}(\F(S))$. This correspondence, which we denote by $P\ra \delta(P)$, is (uniquely) recursively defined by
\begin{equation}
	\delta(x)=\delta_x\ ;\ \delta([P,Q])=[\delta(P),\delta(Q)]\ .
\end{equation}
For $C=\sum_{n\geq 0}C_n\in \ls{\C}{X}$ with $C_n\in \lp{\C}{X}_{n}$, we remark that the sequence $\scal{S\sum_{0\leq n\leq N}C_n}{w}$ is stable (for large $N$). Set
$\delta_{\leq N}:=\delta(\sum_{0\leq n\leq N}C_n)$. We see that $\delta_{\leq N}$ is stable (for large $N$) on every $\F_\al$; we call its limit  $\delta(C)$. It is clear that this limit is a derivation and that it corresponds to $C$.\\
$iv)\Longrightarrow v)$ For every $C=\sum_{n\geq 0}C_n\in \ls{\C}{X}$, the exponential $e^C$ defines a mapping $\phi\in \End(\F(S))$ as indeed $e^{\delta_{\leq N}}$ is stationnary. It is easily checked that this mapping is an automorphism of algebra of $\F(S)$.\\
$v)\Longrightarrow iii)$ For $C_i\in \ls{\C}{X};i=1,2$ which commute we have
\begin{equation}
Se^{C_1}e^{C_2}=\phi_{C_1}(S)e^{C_2}=\phi_{C_1}(Se^{C_2})=\phi_{C_1}\phi_{C_2}(S).
\end{equation}
This proves the existence, for a $C\in \ls{\C}{X}$, of a one-parameter (rational) group $\phi_C^{t}$ in $\Aut(\F(S))$ such that $Se^{tC}=\phi_C^{t}(S)$. This one-parameter (rational) group can be extended to $\R$ as continuity is easily checked by taking the scalar products
$\scal{\phi_C^{t}(S)}{w}=\scal{Se^{tC}}{w}$ and it suffices to specialize the result to $C=x$.\\
$iii)\Longrightarrow ii)$ By stationary limits one has
\begin{equation}
	\scal{Sx}{w}=\lim_{t\ra 0}\frac{1}{t}(\scal{Se^{tx}}{w}-\scal{S}{w})=
	\lim_{t\ra 0}\frac{1}{t}(\scal{\phi_x^{t}(S)}{w}-\scal{S}{w})
\end{equation}
$v)\Longrightarrow i)$ Let $x\in X, t\in \R$, we take $C=tx$ and $\phi_t\in \Aut(\F(S))$  s.t. $\phi_t(S)=Se^{tx}$. It there is $P\in \ncp{\C}{X}$ such that $\scal{S}{P}=0$ one has
\begin{equation}
0=\scal{S}{P}=\phi_t(\scal{S}{P})=\scal{\phi_t(S)}{P}=\scal{Se^{tx}}{P}=\sum_{n=0}^{\deg(P)}\frac{t^n}{n!}\scal{Sx^n}{P}
\end{equation}
and then, for all $z\in V$, the polynomial
\begin{equation}
\sum_{n=0}^{\deg(P)}\frac{t^n}{n!}\scal{S(z)x^n}{P}
\end{equation}
is identically zero over $\R$ hence so are all of its coeefficients in particular $\scal{S(z)x}{P}$ for all $z\in V$. This proves the claim.\\
$i)\Longrightarrow vi)$ Let $P\in ker_\C(S)$ if $P\not=0$ take it of minimal degree with this property. For all $x\in X$, one has $P\in ker_\C(Sx)$ which means $\scal{Sx}{P}=0$ and then $Px^\dagger=0$ as $\deg(Px^\dagger)=\deg(P)-1$. The reconstruction lemma implies that
\begin{equation}
	P=\scal{P}{1}+\sum_{x\in X} (Px^\dagger)x=\scal{P}{1}
\end{equation}
Then, one has $0=\scal{S}{P}=\scal{S}{1}\scal{P}{1}=\scal{P}{1}$ which shows that $ker_\C(S)=\{0\}$. This is equivalent to the statement (vi). \\
$vi)\Longrightarrow i)$ Is obvious as $ker_\C(S)=\{0\}$.\\
\cqfd

\begin{remark} The derivations $\delta_x$ cannot in general be expressed as restrictions of 
derivations of $\H$. 
For example, with equation \mref{diff_eq3}, one has $\delta_{x_0}(\frac{\log(z)^{n+1}}{(n+1)!})=\frac{\log(z)^{n}}{n!}$ but 
$\delta_{x_0}(\scal{S}{ux_1})=0$.
\end{remark}

\section{Conclusion}

In this paper we showed that by using fields of germs, some difficult results can be considerably simplified and extended. For instance, polylogarithms were known to be independant over either $\C[z,1/z,1/(1-z)]$ or, presumably, over "functions which do not involve monodromy"; these two results are now encompassed by Theorem \mref{ind_lin}. We believe that this procedure is not only of theoretical importance, but can be taken into account at the  computational level because every formula (especially analytic) carries with it its domain of validity. As a matter of fact, having at hand the linear independence of coordinate functions over large rings allows one to express uniquely solutions of systems like \mref{nonlinear} in the basis of hyperlogarithms.\\
A valuable  prospect  would be to determine the asymptotic expansion at infinity of the Taylor coefficients of the $y(z)$ as given in \mref{output} for the general case. This  has  been done already   for  the case of singularities at $\{0,1\}$ and for different purposes (see arXiv:1011.0523v2 and http://fr.arxiv.org/abs/0910.1932).\\
\newpage


\begin{thebibliography}{ABC}
%
\bibitem{lappo}{J.A.~Lappo-Danilevsky}.--
	{\em  Th\'eorie des syst\`emes des \'equations diff\'erentielles lin\'eaires},
 Chelsea, New York, 1953.
%
\bibitem{lewin1}{L.~Lewin}.--
	{\em Polylogarithms and associated functions},
 North Holland, New York and Oxford, 1981.
%
\bibitem{dyson}{F.J.~Dyson}.--
 {\em The radiation theories of Tomonaga, Schwinger and Feynman},
 Physical Rev, vol 75, 1949, pp. 486-502.
%
\bibitem{magnus}{W.~Magnus}.--
 On the Exponential Solution of Differential Equation for a Linear Operator,
 Comm. Pure Appl. Math. 7, 1954, pp. 649-673.
%
\bibitem{chen}{K.T.~Chen}.--
	{\em Iterated path integrals},
 Bull. Amer. Math.\allowbreak
 Soc., vol 83, 1977,  pp. 831-879.
%
\bibitem{fliess1}{M.~Fliess}.--
	{\em Fonctionnelles causales non lin\'eaires et ind\'e\-ter\-mi\-n\'ees non commutatives},
 Bull. Soc. Math. France, N$^{\circ}$109, 1981, pp. 3-40.
%
\bibitem{fliess2}{M.~Fliess}.--
	{\em R\'ealisation locale des sys\-t\`e\-mes non li\-n\'e\-aires, al\-g\`e\-bres de Lie fil\-tr\'ees
 transitives et s\'e\-ries g\'en\'e\-ratrices},
 Invent. Math., t 71, 1983, pp. 521-537.
%
\bibitem{feynman}{R.~Feynman}.--
	{\em An Operator Calculus Having Applications in Quantum Electrodynamics},
 Phys. Rev. 84, 108–128 (1951).
%
\bibitem{singervanderput}{M.~van der Put and M.F.~Singer}.--
Galois Theory of Linear Differential Equation, Comprehensive Studies
in Mathematics, vol. 328, Springer-Verlag, Berlin, (2003).
%
\bibitem{wechsung}{G.~Wechsung}.--
Functional Equations of Hyperlogarithms, in \cite{lewin2}
%
\bibitem{FPSAC98}{Hoang Ngoc Minh, M. Petitot and  J. Van der Hoeven}.--
   Polylogarithms and Shuffle Algebra,
   \textit{Proceedings of FPSAC'98}, 1998.
%
\bibitem{BR} J. Berstel, C. Reutenauer, {\it Rational series and their
languages}. EATCS Monographs on Theoretical Computer Science, Springer, 1988.
%
\bibitem{D21} {\sc Duchamp G., Reutenauer C.}, {\it Un crit\`ere de rationalit\'e
provenant de la g\'eom\'etrie noncommutative} Invent. Math. {\bf 128} 613-622.
(1997).
%
\bibitem{Reu} C. Reutenauer, {\it Free Lie algebras}, Oxford University Press, 1993.
%
\bibitem{B_Alg_III} N. Bourbaki, {\it Algebra, chapter III}, Springer (1970)
%
\bibitem{B_Set} N. Bourbaki, {\it Theory of sets}, Springer (2004).
%
\bibitem{lewin2}{L.~Lewin}.--
	{\em  Structural properties of polylogarithms}, Mathematical survey and monographs,
 Amer. Math. Soc., vol 37, 1992.
%
\end{thebibliography}
\end{document}